January 15, 2015

# Power Series with Binomial Sums and Asymptotic Expansions

Khristo N. Boyadzhiev

Ohio Northern University, Department of Mathematics and Statistics

Ada, OH 45810, USA

k-boyadzhiev@onu.edu

**Abstract**

This paper is a study of power series, where the coefficients are binomial expressions (iterated finite differences). Our results can be used for series summation, for series transformation, or for asymptotic expansions involving Stirling numbers of the second kind. In certain cases we obtain asymptotic expansions involving Bernoulli polynomials, poly-Bernoulli polynomials, or Euler polynomials. We also discuss connections to Euler series transformations and other series transformation formulas.

**1 Introduction**

We develop the theory of power series where the coefficients are binomial expressions. In particular, we extend some known results involving classical special functions and polynomials. Guillera and Sondow [9] discussed several interesting representations involving binomial transforms. These include the formulas

$$(1-s)\zeta(s,a) = \sum_{n=0}^{\infty} \frac{1}{n+1} \sum_{k=0}^{n} \binom{n}{k} \frac{(-1)^k}{(k+a)^{s-1}} \qquad (1.1)$$

$$B_m(a) = \sum_{n=0}^{m} \frac{1}{n+1} \sum_{k=0}^{n} \binom{n}{k} (-1)^k (k+a)^m \qquad (1.2)$$

$$\eta(s) = \sum_{n=0}^{\infty} \frac{1}{2^{n+1}} \sum_{k=0}^{n} \binom{n}{k} \frac{(-1)^k}{(k+1)^s} \qquad (1.3)$$

$$E_m(a) = \sum_{n=0}^{\infty} \frac{1}{2^n} \sum_{k=0}^{n} \binom{n}{k} (-1)^k (k+a)^m \qquad (1.4)$$

and some others, where

$$\eta(s) = \sum_{n=1}^{\infty} \frac{(-1)^{n-1}}{n^s} \quad \text{and} \quad \zeta(s,a) = \sum_{n=0}^{\infty} \frac{1}{(n+a)^s} \qquad (1.5)$$

are the Euler eta function and the Hurwitz zeta functions correspondingly, and $B_m(a)$, $E_m(a)$ are the Bernoulli and Euler polynomials. Interesting work in this direction was recorded also by Connon [7]. The representations (1.1) and (1.3) originate from Hasse [11]. The representation (1.2) results from (1.1), since $m\zeta(1-m,a) = -B_m(a)$ for $m = 1, 2, \ldots$. Independently, (1.2) was discovered by Todorov in his detailed investigation of Bernoulli polynomials [16]. This representation is extended in Section 2 to poly-Bernoulli polynomials (see below (2.23)).

Here we want to look at these formulas from a different perspective, including them in a larger theory and connecting them to the series transformation formulas of Euler [2] and the series transformation formulas considered in [3]. Our main results are organized in one proposition, three theorems, and several corollaries. As we shall see, the formalism that we develop leads naturally to certain asymptotic expansions. In section 5 we present new proofs of some classical asymptotic expansions. We also obtain there the asymptotic expansion of the recently introduced Arakawa-Kaneko zeta function (Example 7 in Section 5). Most remarkably, the poly-Bernoulli polynomials and the Bernoulli polynomials appear in several important asymptotic expansions for a large class of functions – see Corollary 2 in Section 2.

## 2 Main results

We study power series whose coefficients are binomial tsums of the form

$$\sum_{k=0}^{n}\binom{n}{k}(-1)^{k} f(y+zk). \tag{2.1}$$

Here $y, z$ are parameters and $f(t)$ is a formal power series

$$f(t) = a_0 + a_1 t + a_2 t^2 + \ldots \tag{2.2}$$

Recall now that for any $m, n$ nonnegative integers, the numbers

$$S(m,n) = \frac{(-1)^n}{n!} \sum_{k=0}^{n}\binom{n}{k}(-1)^{k} k^m$$

are the Stirling numbers of the second kind. Good references for these numbers are the books of Comtet [6], Graham et all [10], and Jordan [12].

**Lemma1**. *For any power series $f(t)$ as in (2.2) and for any integer $n \geq 0$ we have*

$$\sum_{k=0}^{n}\binom{n}{k}(-1)^{k} f(y+zk) = (-1)^n n! \sum_{m=0}^{\infty} a_m \left\{ \sum_{p=0}^{m} S(p,n) z^p y^{m-p} \right\}, \tag{2.3}$$

*where $y, z$ are parameters. In particular, for $y = 0$,*

$$\sum_{k=0}^{n}\binom{n}{k}(-1)^{k} f(zk) = (-1)^n n! \sum_{m=0}^{\infty} a_m z^m S(m,n). \tag{2.4}$$

*Proof.* Starting from the binomial expansion

$$(y+zk)^m = \sum_{p=0}^{m}\binom{m}{p} z^p k^p y^{m-p}$$

we have by changing the order of summation

$$\sum_{k=0}^{n}\binom{n}{k}(-1)^k(y+zk)^m = \sum_{k=0}^{n}\binom{n}{k}(-1)^k\left\{\sum_{p=0}^{m}\binom{m}{p}z^p k^p y^{m-p}\right\}$$

$$= \sum_{p=0}^{m}\binom{m}{p}z^p y^{m-p}\left\{\sum_{k=0}^{n}\binom{n}{k}(-1)^k k^p\right\} = (-1)^n n!\sum_{p=0}^{m}\binom{m}{p}z^p y^{m-p} S(p,n).$$

That is,

$$\sum_{k=0}^{n}\binom{n}{k}(-1)^k(y+zk)^m = (-1)^n n!\sum_{p=0}^{m}\binom{m}{p}S(p,n)z^p y^{m-p}. \tag{2.5}$$

Multiplying both sides by $a_m$ and summing for $m=0,1,2...$ we obtain (2.3). □

**Remark 1**. Note that the summation in the series on the RHS in (2.3) and (2,4) starts, in fact, from $m=n$, since $S(m,n)=0$ when $m<n$. It follows from (2.5) that for any $y,z$,

$$\sum_{k=0}^{n}\binom{n}{k}(-1)^k(y+zk)^m = 0 \quad (m<n),$$

and the sum equals $(-1)^n n! z^n$ for $m=n$.

When $f(t)$ is analytic in some disk centered at the origin with radius $r>0$, the binomial sum on the LHS of (2.4) is analytic in the disk with radius $r/n$ (when $n>0$) and the power series in $z$ on the RHS is absolutely convergent in that disk, being the Taylor series of the LHS.

**Remark 2.** We want to point out that the sums (2.1) are iterated finite differences. Namely, if $\Delta f(t) = f(t+1) - f(t)$, then

$$\sum_{k=0}^{n}\binom{n}{k}(-1)^k f(y+k) = (-1)^n \Delta^n f(y).$$

Recalling that the Newton series of $f(z)$ has the form

$$f(z) = \sum_{k=0}^{\infty}\frac{\Delta^k f(0)}{k!}z(z-1)(z-2)...(z-k+1),$$

we conclude that the above lemma, among other things, gives a connection between the Newton series coefficients of $f(z)$ and its Taylor coefficients. This relationship was observed by Jordan in [12, pp.189-190]. In our opinion formula (2.3) is very important and needs further study and applications.

Lemma 1 is essential for our results. Its first application is in Proposition 1 below, which can be used for series evaluation. For this proposition we need also another lemma.

**Lemma 2**. *For every integer $m > 1$ we have*

$$\sum_{n=1}^{m} S(m,n)(n-1)!\,(-1)^n = 0, \tag{2.6}$$

*and the value of this sum is $-1$ when $m = 1$.*

*Proof.* The Stirling numbers of the second kind satisfy the recurrence relation

$$S(m,n) = n\,S(m-1,n) + S(m-1,n-1),\ 1 \le n \le m$$

and the lemma follows from here by a simple computation. Let $m \ge 2$. Then

$$\sum_{n=1}^{m} S(m,n)(n-1)!\,(-1)^n = \sum_{n=1}^{m} S(m-1,n)\,n!\,(-1)^n + \sum_{n=1}^{m} S(m-1,n-1)(n-1)!\,(-1)^n$$

$$= \sum_{n=1}^{m-1} S(m-1,n)\,n!\,(-1)^n + \sum_{k=0}^{m-1} S(m-1,k)\,k!\,(-1)^{k+1}$$

(by setting $n-1 = k$ and noticing that $S(m-1,m) = 0$)

$$= \sum_{n=1}^{m-1} S(m-1,n)\,n!\,(-1)^n - \sum_{k=1}^{m-1} S(m-1,k)\,k!\,(-1)^k - S(m-1,0) = 0$$

The case $m = 1$ is obvious and the proof is completed. □

From the above two lemmas we derive the proposition.

**Proposition 1**. *Let $f(t)$ be a formal power series as in* (2.2). *Then for every $z$ we have*

$$\sum_{n=1}^{\infty}\frac{1}{n}\left\{\sum_{k=0}^{n}\binom{n}{k}(-1)^k f(zk)\right\} = -f'(0)z \qquad (2.7)$$

*Proof.* We multiply both sides in equation (2.4) by $\frac{1}{n}$ and sum for $n = 1, 2, ...$ . This yields

$$\sum_{n=1}^{\infty}\frac{1}{n}\left\{\sum_{k=0}^{n}\binom{n}{k}(-1)^k f(zk)\right\} = \sum_{m=0}^{\infty}a_m z^m\left\{\sum_{n=1}^{m}S(m,n)(n-1)!(-1)^n\right\} = -a_1 z,$$

which is the desired result. (see also (2.19) below). □

When $f(t)$ is a polynomial of degree $p$, the series on the LHS in (2.7) truncates, as the binomial sums become zeros for $n > p$ (see Remark 1). A large class of functions for which (2.7) holds is presented in [4] together with several applications.

We continue with results related to lemma 1.

**Example 1.** Consider the polynomial

$$g(t) = \binom{t}{p} = \frac{1}{p!}\sum_{m=0}^{p}s(p,m)t^m, \qquad (2.8)$$

where $p \geq 0$ is an integer and $s(m, p)$ are the Stirling numbers of the first kind. From Lemma 1,

$$\sum_{k=0}^{n}\binom{n}{k}(-1)^k\binom{zk}{p} = (-1)^n\frac{n!}{p!}\sum_{m=0}^{p}s(p,m)S(m,n)z^m. \qquad (2.9)$$

These sums appeared in Todorov's paper [17], in the Taylor series expansion

$$\left((1+t)^z - 1\right)^n = \sum_{p=0}^{\infty}t^p\left\{(-1)^n\sum_{k=0}^{n}\binom{n}{k}(-1)^k\binom{zk}{p}\right\}.$$

When $p < n$, both sides in (2.9) are zeros.

Now consider the exponential polynomials [3], [5]

$$\varphi_m(x) = \sum_{n=0}^{m} S(m,n) x^n$$

In the following theorem the notation BE stands for Binomial Exponential Series.

**Theorem 1.** *With $f(t)$ as in (2.2) and with parameters $y$, $z$, the following representation holds*

$$\sum_{n=0}^{\infty} \frac{x^n}{n!} \left\{ \sum_{k=0}^{n} \binom{n}{k} (-1)^k f(y+zk) \right\} = \sum_{m=0}^{\infty} a_m \left\{ \sum_{p=0}^{m} \binom{m}{p} z^p y^{m-p} \varphi_p(-x) \right\}. \qquad (2.10)$$

*When $y = 0$ this becomes*

**BE** $\qquad \sum_{n=0}^{\infty} \frac{x^n}{n!} \left\{ \sum_{k=0}^{n} \binom{n}{k} (-1)^k f(zk) \right\} = \sum_{m=0}^{\infty} a_m \varphi_m(-x) z^m \qquad (2.11)$

*Proof.* Multiplying both sides of (2.3) by $\dfrac{x^n}{n!}$ and summing for $n = 0, 1 \ldots$ we come to (2.10) after changing the order of summation on the RHS. Note that since $S(p, n) = 0$ for $n > p$, the summation on the RHS leads exactly to the exponential polynomials $\varphi_p(-x)$. □

We do not specify convergence in the above formulas. If $f(t)$ is an entire function, the partial sums of the series on the LHS in (2.11) are defined for every $z$. If the series on the RHS converges for some $z$, the formula will hold for all $z$ in some neighborhood of the origin. Depending on the choice of the function $f(t)$, the series on the RHS in the above equations may diverge. As we shall see later, in some cases the RHS represents an asymptotic expansion as explained in [3].

From the above theorem with $x = -1$ we deduce the immediate corollary.

**Corollary 1.** *For any power series $f(t)$ as in (2.2) we have*

$$\sum_{n=0}^{\infty} \frac{(-1)^n}{n!} \left\{ \sum_{k=0}^{n} \binom{n}{k} (-1)^k f(zk) \right\} = \sum_{m=0}^{\infty} a_m b_m z^m, \qquad (2.12)$$

*where the numbers*

$$b_n = \varphi_n(1) = \sum_{k=0}^{n} S(n,k)$$

*are the well-known Bell numbers* [5].

**Example 2.** From (2.9) and Theorem 1 we obtain for every integer $p \geq 0$,

$$\sum_{n=0}^{p} \frac{x^n}{n!} \left\{ \sum_{k=0}^{n} \binom{n}{k}(-1)^k \binom{zk}{p} \right\} = \frac{1}{p!} \sum_{m=0}^{p} s(p,m)\, \varphi_m(-x)\, z^m \, . \tag{2.13}$$

We shall prove now analogous expansions involving another sequence of polynomials. Let

$$\omega_m(x) = \sum_{n=0}^{m} S(m,n) n! x^n ,$$

($m = 0, 1, ...$) be the geometric polynomials discussed in [3]. The geometric polynomials $\omega_m$ will replace the exponential polynomials $\varphi_m$ on the RHS of the new formulas. Namely, multiplying both sides of (2.3) by $x^n$ and summing for $n = 0, 1, ...$, we obtain (after changing the order of summation) the following theorem.

**Theorem 2**. *With $f(t)$ as in* (2.2) *and with $y, z$ parameters,*

$$\sum_{n=0}^{\infty} x^n \left\{ \sum_{k=0}^{n} \binom{n}{k}(-1)^k f(y+zk) \right\} = \sum_{m=0}^{\infty} a_m \left\{ \sum_{p=0}^{m} \binom{m}{p} z^p y^{m-p} \omega_p(-x) \right\} . \tag{2.14}$$

*When $y = 0$,*

**BG** $$\sum_{n=0}^{\infty} x^n \left\{ \sum_{k=0}^{n} \binom{n}{k}(-1)^k f(zk) \right\} = \sum_{m=0}^{\infty} a_m z^m \omega_m(-x) \tag{2.15}$$

(BG stands for Binomial Geometric Series).

**Example 3**. From (2.9) and (2.15) we have for all integers $p \geq 0$,

$$\sum_{n=0}^{p} x^n \left\{ \sum_{k=0}^{n} \binom{n}{k} (-1)^k \binom{zk}{p} \right\} = \frac{1}{p!} \sum_{m=0}^{p} s(p,m)\, \omega_m(-x)\, z^m \ . \qquad (2.16)$$

Next we list some variations of these formulas. Integrating for $x$ in (2.14) yields

$$\sum_{n=0}^{\infty} \frac{x^{n+1}}{n+1} \left\{ \sum_{k=0}^{n} \binom{n}{k} (-1)^k f(y+zk) \right\} = \sum_{m=0}^{\infty} a_m \left\{ \sum_{p=0}^{m} \binom{m}{p} z^p y^{m-p} \int_0^x \omega_p(-t)\, dt \right\}$$

or, explicitly, after evaluating the integral on the RHS and reducing by $x$,

$$\sum_{n=0}^{\infty} \frac{x^n}{n+1} \left\{ \sum_{k=0}^{n} \binom{n}{k} (-1)^k f(y+zk) \right\}$$

$$= \sum_{m=0}^{\infty} a_m \left\{ \sum_{p=0}^{m} \binom{m}{p} z^p y^{m-p} \sum_{j=0}^{p} S(p,j)\, j!\, \frac{(-x)^j}{j+1} \right\}.$$

We can integrate here for $x$, then reduce by x, and integrate again. Repeating this operation $r-1$ times yields for any $r \geq 1$

$$\sum_{n=0}^{\infty} \frac{x^n}{(n+1)^r} \left\{ \sum_{k=0}^{n} \binom{n}{k} (-1)^k f(y+zk) \right\} \qquad (2.17)$$

$$= \sum_{m=0}^{\infty} a_m \left\{ \sum_{p=0}^{m} \binom{m}{p} z^p y^{m-p} \sum_{j=0}^{p} S(p,j)\, j!\, \frac{(-x)^j}{(j+1)^r} \right\}.$$

When $y = 0$ this becomes

$$\sum_{n=0}^{\infty} \frac{x^n}{(n+1)^r} \left\{ \sum_{k=0}^{n} \binom{n}{k} (-1)^k f(zk) \right\} = \sum_{m=0}^{\infty} a_m z^m \left\{ \sum_{j=0}^{m} S(m,j)\, j!\, \frac{(-x)^j}{(j+1)^r} \right\}. \qquad (2.18)$$

Dividing by $x$ in (2.15) and integrating (removing first the terms without $x$, i.e. $a_0$ from both sides) yields

$$\sum_{n=1}^{\infty} \frac{x^n}{n} \left\{ \sum_{k=0}^{n} \binom{n}{k} (-1)^k f(zk) \right\} = \sum_{m=0}^{\infty} a_m z^m \left\{ \sum_{p=1}^{m} S(m,p)(p-1)!(-x)^p \right\}. \qquad (2.19)$$

We want to show now some immediate corollaries involving the Bernoulli numbers $B_m$, the Bernoulli polynomials,

$$B_m(y) = \sum_{j=0}^{m} \binom{m}{j} y^{m-j} B_j,$$

as well as the poly-Bernoulli numbers $B_n^{(q)}$, and the poly-Bernoulli polynomials $B_n^{(q)}(y)$. The poly-Bernoulli numbers were introduced by Kaneko [13] through the generating function

$$\frac{\text{Li}_q(1-e^{-t})}{1-e^{-t}} = \sum_{n=0}^{\infty} B_n^{(q)} \frac{t^n}{n!}, \qquad (2.20)$$

where

$$\text{Li}_q(x) = \sum_{n=1}^{\infty} \frac{x^n}{n^q}$$

is the polylogarithm. Kaneko showed that for $q=1$ we have $B_n^{(1)} = (-1)^n B_n$, i.e. $B_1^{(1)} = 1/2$ and $B_n^{(1)} = B_n$ $(n \neq 1)$. Also proved in [13] was the representation of the poly-Bernoulli numbers in terms of the Stirling numbers of the second kind,

$$B_n^{(q)} = (-1)^n \sum_{j=0}^{n} S(n,j) j! \frac{(-1)^j}{(j+1)^q}.$$

This extends the well known representation of Bernoulli numbers in terms of $S(n,j)$, namely,

$$B_n = \sum_{j=0}^{n} S(n,j) j! \frac{(-1)^j}{j+1}.$$

In analogy to Bernoulli polynomials one can define the poly-Bernoulli polynomials by

$$B_m^{(q)}(y) = \sum_{j=0}^{m} \binom{m}{j} y^{m-j} B_j^{(q)}.$$

This is equivalent to the definition given by Bayad and Hamahata in [1] in terms of the generating function

$$\frac{\text{Li}_q(1-e^{-t})}{1-e^{-t}} e^{xt} = \sum_{n=0}^{\infty} B_n^{(q)}(x) \frac{t^n}{n!}.$$

The polynomials $B_m^{(q)}(x)$ slightly differ from the similar polynomials defined in [8].

When $q = 1$ we have $B_m^{(1)}(y) = (-1)^m B_m(-y)$, since

$$B_m^{(1)}(y) = \sum_{j=0}^{m} \binom{m}{j} y^{m-j} B_j = \sum_{j=0}^{m} \binom{m}{j} y^{m-j} (-1)^{-j} B_j = (-1)^m B_m(-y).$$

We shall use now equation (2.17) to obtain representations involving the Bernoulli and the poly-Bernoulli numbers and polynomials.

**Corollary 2.** *For every power series (2.2), every integer $r \geq 1$, and with parameters $y$ and $z$*

$$\sum_{n=0}^{\infty} \frac{1}{(n+1)^r} \left\{ \sum_{k=0}^{n} \binom{n}{k} (-1)^k f(y+zk) \right\} = \sum_{m=0}^{\infty} a_m (-z)^m B_m^{(r)}(-y z^{-1}). \quad (2.21)$$

*When $y = 0$,*

$$\sum_{n=0}^{\infty} \frac{1}{(n+1)^r} \left\{ \sum_{k=0}^{n} \binom{n}{k} (-1)^k f(zk) \right\} = \sum_{m=0}^{\infty} a_m (-z)^m B_m^{(r)}. \quad (2.22)$$

*In particular, with $r = 1$ we have the following expansions involving the regular Bernoulli numbers and polynomials,*

$$\sum_{n=0}^{\infty} \frac{1}{n+1} \left\{ \sum_{k=0}^{n} \binom{n}{k} (-1)^k f(y+zk) \right\} = \sum_{m=0}^{\infty} a_m z^m B_m(y z^{-1}) \quad (2.23)$$

$$\sum_{n=0}^{\infty}\frac{1}{n+1}\left\{\sum_{k=0}^{n}\binom{n}{k}(-1)^k f(zk)\right\}=\sum_{m=0}^{\infty}a_m z^m B_m. \qquad (2.24)$$

With the choice $f(t)=t^m$ and $z=1$ in (2.23) we obtain the representation (1.2). An analogous representation for the poly-Bernoulli polynomials results from (2.21) with $z=-1$

$$B_m^{(r)}(y)=\sum_{n=0}^{m}\frac{1}{(n+1)^r}\sum_{k=0}^{n}\binom{n}{k}(-1)^k (y-k)^m \qquad (2.25)$$

(proved by a different method in [1]). Note that the summation goes to $n=m$, since for $n>m$

$$\sum_{k=0}^{n}\binom{n}{k}(-1)^k (y-k)^m = 0.$$

For the function $f(t)=e^t$ we have

$$\sum_{k=0}^{n}\binom{n}{k}(-1)^k f(zk)=(1-e^z)^n$$

and the LHS of (2.24) becomes

$$\frac{1}{1-e^z}\sum_{n=0}^{\infty}\frac{(1-e^z)^{n+1}}{n+1}=\frac{1}{1-e^z}\left(-\log(1-(1-e^z))\right)=\frac{-z}{1-e^z}$$

(for convergence we need $\mathrm{Re}\, z<0$). Therefore, (2.24) provides the expected representation

$$\frac{z}{e^z-1}=\sum_{n=0}^{\infty}\frac{z^n}{n!}B_n,$$

which converges for $|z|<2\pi$. This is a good example how we need to adjust sometimes $z$ or $x,y$ in the above theorems in order to assure convergence.

*Proof of the corollary.* With $y,z$ parameters we first compute

$$\sum_{p=0}^{m}\binom{m}{p} z^p y^{m-p} \sum_{j=0}^{p} S(p,j) j! \frac{(-1)^j}{(j+1)^r}$$

$$= (-1)^m z^m \sum_{p=0}^{m}\binom{m}{p}\left(\frac{-y}{z}\right)^{m-p} (-1)^p \sum_{j=0}^{p} S(p,j) j! \frac{(-1)^j}{(j+1)^r}$$

$$= (-1)^m z^m \sum_{p=0}^{m}\binom{m}{p}\left(\frac{-y}{z}\right)^{m-p} B_p^{(r)} = (-z)^m B_m^{(r)}(-y z^{-1}).$$

That is, we have the equation

$$\sum_{p=0}^{m}\binom{m}{p} z^p y^{m-p} \sum_{j=0}^{p} S(p,j) j! \frac{(-1)^j}{(j+1)^r} = (-z)^m B_m^{(r)}(-y z^{-1}).$$

Multiplying this equation by $a_m$ and summing for $m = 0, 1, \ldots$ we obtain

$$\sum_{m=0}^{\infty} a_m \left\{ \sum_{p=0}^{m}\binom{m}{p} z^p y^{m-p} \sum_{j=0}^{p} S(p,j) j! \frac{(-1)^j}{(j+1)^r} \right\} = \sum_{m=0}^{\infty} a_m (-z)^m B_m^{(r)}(-y z^{-1})$$

and (2.21) follows from here and (2.17). □

Notice that with the selection $f(t) = 1 - e^{-t}$ in (2.22) we arrive at equation (2.20).

It is known that

$$\omega_m\left(\frac{-1}{2}\right) = \sum_{p=0}^{m}(-1)^p \frac{p!}{2^p} = \frac{2}{m+1}\left(1 - 2^{m+1}\right) B_{m+1} \tag{2.26}$$

(see, for instance [10, Problem 6.76, page 317]). Setting $x = \frac{1}{2}$ in (2.15) we obtain one more interesting representation.

**Corollary 3.** *With $f(t)$ as in* (2.2),

$$\sum_{n=0}^{\infty}\frac{1}{2^{n+1}}\left\{\sum_{k=0}^{n}\binom{n}{k}(-1)^k f(zk)\right\}=\sum_{m=0}^{\infty}\frac{1-2^{m+1}}{m+1}a_m B_{m+1} z^m. \quad (2.27)$$

Again, we do not specify convergence in the above series. When $f(t)$ is a polynomial, all formulas in Theorems 1 and 2, as well as the formulas in Corollaries 2 and 3 provide closed form evaluations of the series on the LHS. By choosing different functions $f(t)$ in Corollaries 2 and 3 we can generate various asymptotic expansions. It is really amazing that the poly-Bernoulli and the Bernoulli polynomials are present in all these situations.

## 3 Series with Euler polynomials

We shall record here another corollary involving this time Euler's numbers $E_m$ and Euler's polynomials $E_m(x)$. First we need to show a special connection between the geometric polynomials and Euler's polynomials. The geometric polynomials $\omega_m(t)$ are related to the geometric series in the following way: For every $|x|<1$ and every $m=0,1,2...$, we have

$$\sum_{n=0}^{\infty} n^m x^n = \frac{1}{1-x}\omega_m\left(\frac{x}{1-x}\right) \quad (3.1)$$

(see [3]). We shall use this property to find the Taylor series of the function

$$h(t)=\frac{1}{\mu e^{\lambda t}+1},$$

where $\lambda, \mu$ are parameters. For this purpose we need to evaluate the higher derivatives of $h$ at zero. When $|\mu e^{\lambda t}|<1$ we expand as geometric series

$$\frac{1}{\mu e^{\lambda t}+1}=\frac{1}{1-(-\mu e^{\lambda t})}=\sum_{n=0}^{\infty}(-\mu)^n e^{\lambda t n}.$$

From this

$$\left(\frac{d}{dt}\right)^m \frac{1}{\mu e^{\lambda t}+1}=\lambda^m \sum_{n=0}^{\infty}(-\mu)^n n^m e^{\lambda t n}$$

and in view of (3.1),

$$\left(\frac{d}{dt}\right)^m \frac{1}{\mu e^{\lambda t}+1} = \frac{1}{\mu e^{\lambda t}+1}\omega_m\left(\frac{-\mu e^{\lambda t}}{\mu e^{\lambda t}+1}\right) \tag{3.2}$$

so that

$$h^{(m)}(0) = \left(\frac{d}{dt}\right)^m \frac{1}{\mu e^{\lambda t}+1}\bigg|_{t=0} = \frac{\lambda^m}{\mu+1}\omega_m\left(\frac{-\mu}{\mu+1}\right), \tag{3.4}$$

which provides the Taylor series representation

$$\frac{1}{\mu e^{\lambda t}+1} = \frac{1}{\mu+1}\sum_{m=0}^{\infty}\lambda^m \omega_m\left(\frac{-\mu}{\mu+1}\right)\frac{t^m}{m!}. \tag{3.5}$$

In particular, with $\lambda = \mu = 1$,

$$\frac{2}{e^t+1} = \sum_{m=0}^{\infty}\omega_m\left(\frac{-1}{2}\right)\frac{t^m}{m!}. \tag{3.6}$$

Multiplying both sides by $e^{xt}$ and using Cauchy's rule for multiplication of power series (see also (4.4) below) we can write

$$\frac{2e^{xt}}{e^t+1} = \sum_{m=0}^{\infty}\frac{t^m}{m!}\left\{\sum_{k=0}^{m}\binom{m}{k}\omega_k\left(\frac{-1}{2}\right)x^{m-k}\right\}. \tag{3.7}$$

At the same time, the function on the LHS here is the generating function for Euler's polynomials $E_m(x)$, i.e.

$$\frac{2e^{xt}}{e^t+1} = \sum_{m=0}^{\infty}E_m(x)\frac{t^m}{m!}. \tag{3.8}$$

Comparing both series we find

$$E_m(x) = \sum_{k=0}^{m}\binom{m}{k}\omega_k\left(\frac{-1}{2}\right)x^{m-k}. \tag{3.9}$$

In particular, when $x = 0$,

$$E_m(0) = \omega_m\left(\frac{-1}{2}\right) = \frac{2}{m+1}\left(1 - 2^{m+1}\right)B_{m+1} \tag{3.10}$$

Euler's numbers $E_m$ are defined by

$$E_m = 2^m E_m\left(\frac{1}{2}\right). \tag{3.11}$$

Thus we arrive at the following corollary.

**Corollary 4.** *With $f(t)$ as in (2.2),*

$$\sum_{n=0}^{\infty}\frac{1}{2^n}\left\{\sum_{k=0}^{n}\binom{n}{k}(-1)^k f(y+k)\right\} = \sum_{m=0}^{\infty} a_m E_m(y), \tag{3.12}$$

*where $E_m(y)$ are Euler's polynomial. In particular, with $y = \frac{1}{2}$,*

$$\sum_{n=0}^{\infty}\frac{1}{2^n}\left\{\sum_{k=0}^{n}\binom{n}{k}(-1)^k f\left(\frac{1}{2}+k\right)\right\} = \sum_{m=0}^{\infty}\frac{1}{2^m} a_m E_m. \tag{3.13}$$

*Proof.* In (2,14) we set $x = \frac{1}{2}$ and $z = 1$. Then we use (3.9) for (3.12) and (3.11) for (3.13). □

Applying the above corollary to the function $f(t) = t^m$ yields (1.4).

## 4 Connection to other series transformation formulas

In [3] the present author considered two series transformation formulas (STF). Namely, the Exponential STF,

**ES** $$\sum_{n=0}^{\infty}\frac{x^n}{n!} f(zk) = e^x \sum_{m=0}^{\infty} a_m z^m \varphi_m(x) \tag{4.1}$$

and the Geometric STF (for $|x| < 1$).

**GS** $$\sum_{n=0}^{\infty} x^n f(zk) = \frac{1}{1-x} \sum_{m=0}^{\infty} a_m z^m \omega_m\left(\frac{x}{1-x}\right), \tag{4.2}$$

where $\varphi_m(x)$ and $\omega_m(x)$ are the exponential polynomials and the geometric polynomials correspondingly. It was shown that in some cases the **ES** and **GS** can be used to evaluate power series in a closed form. In other cases these representations provide asymptotic expansions. To show the connection of **ES** and **GS** with the expansions in the preceding section we shall use Euler's series transformations as described in [2]. Let $f(t)$ be defined by (2.2). The **Exponential Euler Series Transformation (EE)** says that

**EE** $$e^x \sum_{n=0}^{\infty} \frac{(-x)^n}{n!} f(zk) = \sum_{n=0}^{\infty} \frac{x^n}{n!} \left\{ \sum_{k=0}^{n} \binom{n}{k} (-1)^k f(zk) \right\} \tag{4.3}$$

or, more generally,

$$e^{\lambda x} \sum_{n=0}^{\infty} \frac{(-x)^n}{n!} f(zk) = \sum_{n=0}^{\infty} \frac{x^n}{n!} \left\{ \sum_{k=0}^{n} \binom{n}{k} (-1)^k \lambda^{n-k} f(zk) \right\}, \tag{4.4}$$

( with $z, \lambda$ parameters). The **geometric** version is ([1], [3])

**EG** $$\frac{1}{1-t} \sum_{n=0}^{\infty} (-1)^n f(zk) \left(\frac{t}{1-t}\right)^n = \sum_{n=0}^{\infty} t^n \left\{ \sum_{k=0}^{n} \binom{n}{k} (-1)^k f(zk) \right\}. \tag{4.5}$$

**Theorem 3.** *Consider the three formulas* **BE**, **ES**, *and* **EE** *(that is, (2.11), (4.1), and (4.3)). Then any two of them together imply the third one. The same is true for* **BG**, **GS**, *and* **EG**.

*Proof.* The proof of the first part is straightforward. In fact, we can unite the formulas together.

$$\sum_{n=0}^{\infty} \frac{x^n}{n!} \left\{ \sum_{k=0}^{n} \binom{n}{k} (-1)^k f(zk) \right\} = e^x \sum_{n=0}^{\infty} \frac{(-x)^n}{n!} f(zk) = \sum_{m=0}^{\infty} a_m z^m \varphi_m(-x). \tag{4.6}$$

For the second part of the theorem, we first use the substitution $x = \frac{-t}{1-t}, t = \frac{-x}{1-x}$, in **GS** to put it in the form

$$\frac{1}{1-t}\sum_{n=0}^{\infty}\left(\frac{-t}{1-t}\right)^{n}f(zk)=\sum_{m=0}^{\infty}a_{m}z^{m}\omega_{m}(-t), \tag{4.7}$$

and then we write (using **EG** for the first equality and (3.5) for the second)

$$\sum_{n=0}^{\infty}t^{n}\left\{\sum_{k=0}^{n}\binom{n}{k}(-1)^{k}f(zk)\right\}=\frac{1}{1-t}\sum_{n=0}^{\infty}(-1)^{n}f(zk)\left(\frac{t}{1-t}\right)^{n}=\sum_{m=0}^{\infty}a_{m}z^{m}\omega_{m}(-t). \tag{4.8}$$

The proof is complete. □

**Example 4.** For any complex number $s$ consider the function

$$f_{s}(t)=\frac{1}{(1+t)^{s}}=\sum_{m=0}^{\infty}\binom{-s}{m}t^{m}, \quad |t|<1.$$

According to (4.6) we have

$$\sum_{n=0}^{\infty}\frac{x^{n}}{n!}\left\{\sum_{k=0}^{n}\binom{n}{k}\frac{(-1)^{k}}{(1+kz)^{s}}\right\}=e^{x}\sum_{n=0}^{\infty}\frac{(-x)^{n}}{n!(1+nz)^{s}}=\sum_{m=0}^{\infty}\binom{-s}{m}z^{m}\varphi_{m}(-x). \tag{4.9}$$

Replacing here $z$ by $\frac{1}{\lambda}$ and $x$ by $-x$ we come to the representation

$$e_{s}(x,\lambda)\equiv\sum_{n=0}^{\infty}\frac{x^{n}}{n!(n+\lambda)^{s}}=\sum_{m=0}^{\infty}\binom{-s}{m}\frac{(-1)^{m}}{\lambda^{m+s}}\varphi_{m}(x), \tag{4.10}$$

which is the asymptotic expansion of the polyexponential function $e_{s}(x,\lambda)$ in terms of $\lambda$. This expansion was obtained in [3].

## 5 Asymptotic expansions

In this section we give new proofs of some classical asymptotic expansions and also provide some new results. In particular, in Example 7 we find the asymptotic expansion of the recently introduced Arakawa-Kaneko zeta function.

**Example 5.** Expansions related to the Lerch Transcendent.

The Lerch Transcendent is defined by the series

$$\Phi(z,s,a) = \sum_{k=0}^{\infty} \frac{z^k}{(k+a)^s}$$

($a > 0$). When $z = 1$ and $\operatorname{Re} s > 1$ this is the Hurwitz zeta function

$$\zeta(s,a) = \Phi(1,s,a) = \sum_{k=0}^{\infty} \frac{1}{(k+a)^s} \,.$$

We have the representations

$$\Phi(x,s,a) = \frac{1}{\Gamma(s)} \int_0^{\infty} \frac{t^{s-1} e^{-at}}{1 - xe^{-t}} dt, \qquad \zeta(s,a) = \frac{1}{\Gamma(s)} \int_0^{\infty} \frac{t^{s-1} e^{-at}}{1 - e^{-t}} dt, \tag{5.1}$$

$$\sum_{k=0}^{n} \binom{n}{k} \frac{(-1)^k x^k}{(k+a)^s} = \frac{1}{\Gamma(s)} \int_0^{\infty} t^{s-1} e^{-at} (1 - xe^{-t})^n \, dt \tag{5.2}$$

$$= \frac{1}{\Gamma(s)} \int_0^{\infty} (1 - xe^{-t})^{n+1} \frac{t^{s-1} e^{-at}}{1 - xe^{-t}} dt \,,$$

and for $0 < x \leq 1$ also,

$$s\,\Phi(x,s+1,a) - \log x\, \Phi(x,s,a) = \sum_{n=0}^{\infty} \frac{1}{n+1} \left\{ \sum_{k=0}^{n} \binom{n}{k} (-1)^k \frac{x^k}{(a+k)^s} \right\}. \tag{5.3}$$

With $x = 1$ this becomes

$$s\,\zeta(s+1,a) = \sum_{n=0}^{\infty} \frac{1}{n+1} \left\{ \sum_{k=0}^{n} \binom{n}{k} \frac{(-1)^k}{(a+k)^s} \right\}, \tag{5.4}$$

also,

$$\Phi(-x, s, a) = \sum_{n=0}^{\infty} \frac{1}{2^{n+1}} \left\{ \sum_{k=0}^{n} \binom{n}{k} (-1)^k \frac{x^k}{(k+a)^s} \right\}, \tag{5.5}$$

and in particular, when $x = 1$,

$$\eta(s,a) = \Phi(-1,s,a) = \sum_{n=0}^{\infty} \frac{1}{2^{n+1}} \left\{ \sum_{k=0}^{n} \binom{n}{k} \frac{(-1)^k}{(k+a)^s} \right\} . \tag{5.6}$$

*Proof.* Starting from the well-known representation

$$\frac{1}{(k+a)^s} = \frac{1}{\Gamma(s)} \int_0^{\infty} t^{s-1} e^{-kt} e^{-at} dt \tag{5.7}$$

we multiply both sides by $x^k$ and sum for $k = 0, 1, ..., \infty$ to get (5.1). In order to prove (5.2) we multiply (5.7) by $\binom{n}{k}(-1)^k x^k$ and sum for $k = 0, 1, ..., n$. From (5.2) and with $0 < x \leq 1$ we have

$$\sum_{n=0}^{\infty} \frac{1}{n+1} \left\{ \sum_{k=0}^{n} \binom{n}{k} (-1)^k \frac{x^k}{(k+a)^s} \right\} = \frac{1}{\Gamma(s)} \int_0^{\infty} \left\{ \sum_{n=0}^{\infty} \frac{(1-xe^{-t})^{n+1}}{n+1} \right\} \frac{t^{s-1} e^{-at}}{1-xe^{-t}} dt$$

$$= \frac{1}{\Gamma(s)} \int_0^{\infty} \{-\log(xe^{-t})\} \frac{t^{s-1} e^{-at}}{1-xe^{-t}} dt = \frac{1}{\Gamma(s)} \int_0^{\infty} \{-\log x + t\} \frac{t^{s-1} e^{-at}}{1-xe^{-t}} dt$$

$$= \frac{-\log x}{\Gamma(s)} \int_0^{\infty} \frac{t^{s-1} e^{-at}}{1-xe^{-t}} dt + \frac{1}{\Gamma(s)} \int_0^{\infty} \frac{t^s e^{-at}}{1-xe^{-t}} dt$$

$$= -\log x \, \Phi(x,s,a) + s \, \Phi(x,s+1,a),$$

which proves (5.3). Next, from (5.2) again,

$$\sum_{n=0}^{\infty} \frac{1}{2^{n+1}} \left\{ \sum_{k=0}^{n} \binom{n}{k} (-1)^k \frac{x^k}{(k+a)^s} \right\} = \frac{1}{\Gamma(s)} \int_0^{\infty} \left\{ \frac{1}{2} \sum_{n=0}^{\infty} \left( \frac{1-xe^{-t}}{2} \right)^n \right\} t^{s-1} e^{-at} dt$$

$$= \frac{1}{\Gamma(s)} \int_0^{\infty} \left\{ \frac{1}{2} \frac{1}{1 - \left( \frac{1-xe^{-t}}{2} \right)} \right\} t^{s-1} e^{-at} dt = \frac{1}{\Gamma(s)} \int_0^{\infty} \left\{ \frac{1}{1+xe^{-t}} \right\} t^{s-1} e^{-at} dt$$

$$= \Phi(-x, s, a),$$

which is (5.5). □

For some similar results see [7].

Now let $0 < x \leq 1$ and consider the function

$$f(t) = \frac{x^t}{(a+t)^s} \qquad (5.8)$$

We write $x^t = e^{t \log x}$ and expand this in power series for $t$. Then using Cauchy's rule for multiplication of power series and setting $m + k = n$ we write

$$\frac{x^t}{(a+t)^s} = \frac{1}{a^s}\left(1 + \frac{t}{a}\right)^{-s} \sum_{k=0}^{\infty} \frac{\log^k x}{k!} t^k = \frac{1}{a^s}\left\{\sum_{m=0}^{\infty} \binom{-s}{m} \frac{t^m}{a^m}\right\}\left\{\sum_{k=0}^{\infty} \frac{\log^k x}{k!} t^k\right\} \qquad (5.9)$$

$$= \sum_{n=0}^{\infty} t^n \left\{\sum_{k=0}^{n} \binom{-s}{n-k} \frac{\log^k x}{k! \, a^{n+s-k}}\right\}$$

When $x = 1$ this becomes

$$f(t) = \frac{1}{(a+t)^s} = \sum_{m=0}^{\infty} \binom{-s}{m} \frac{t^m}{a^{m+s}} \, . \qquad (5.10)$$

We apply now formula (5.3) to obtain (in view of (5.9) and (2.24)) the asymptotic expansion

$$s \, \Phi(x, s+1, a) - \log x \, \Phi(x, s, a) = \sum_{m=0}^{\infty} B_m \left\{\sum_{k=0}^{m} \binom{-s}{m-k} \frac{\log^k x}{k! \, a^{m+s-k}}\right\} \, .$$

In particular, when $x = 1$ only the $k = 0$ term in the second sum is nonzero and we come to the well-known asymptotic series for the Hurwitz zeta function (when $a \to \infty$)

$$s \, \zeta(s+1, a) = \sum_{m=0}^{\infty} \binom{-s}{m} \frac{B_m}{a^{m+s}} \, . \qquad (5.11)$$

(See [14, p. 25] or [15, p. 610]. The expansion is written there in a somewhat different, but equivalent form.)

From (3.12) we also derive the asymptotic expansion

$$\eta(s, y+a) = \sum_{n=0}^{\infty} \frac{1}{2^{n+1}} \left\{ \sum_{k=0}^{n} \binom{n}{k} \frac{(-1)^k}{(k+y+a)^s} \right\} = \frac{1}{2} \sum_{m=0}^{\infty} \binom{-s}{m} \frac{E_m(y)}{a^{m+s}} . \tag{5.12}$$

In particular, when $y = 0$,

$$\eta(s, a) = \frac{1}{2} \sum_{m=0}^{\infty} \binom{-s}{m} \frac{E_m(0)}{a^{m+s}} = \sum_{m=0}^{\infty} \binom{-s}{m} \frac{(1-2^{m+1})B_{m+1}}{(m+1) a^{m+s}} . \tag{5.13}$$

**Example 6.** The log-gamma and digamma functions.

Here we show how the results in Corollary 2 can be used to derive the classical asymptotic expansions of the log-gamma function $\log \Gamma(z)$ and the digamma function $\psi(z) = \frac{d}{dz} \log \Gamma(z)$ (see [18]). We shall use the integral representation

$$\psi(z) = \log z + \int_0^{\infty} \left\{ \frac{1}{x} - \frac{1}{1-e^{-x}} \right\} e^{-zx} dx . \tag{5.14}$$

Let $\operatorname{Re} z > 0$ and consider the function

$$f(t) = \log\left(1 + \frac{t}{z}\right) \tag{5.15}$$

for $t \geq 0$. Using the representation

$$\log\left(1 + \frac{t}{z}\right) = \int_0^{\infty} \frac{1 - e^{-tx}}{x} e^{-zx} dx ,$$

we find

$$\sum_{k=0}^{n} \binom{n}{k} (-1)^k \log\left(1 + \frac{k}{z}\right) = \int_0^{\infty} \left\{ \sum_{k=0}^{n} \binom{n}{k} (-1)^k - \sum_{k=0}^{n} \binom{n}{k} (-1)^k e^{-kx} \right\} \frac{e^{-zx}}{x} dx . \tag{5.16}$$

Note that when $n = 0$ both sides are zeros. For $n \geq 1$ we have

$$\sum_{k=0}^{n} \binom{n}{k} (-1)^k \log\left(1 + \frac{k}{z}\right) = \int_0^{\infty} \left\{ -\sum_{k=0}^{n} \binom{n}{k} (-1)^k e^{-kx} \right\} \frac{e^{-zx}}{x} dx = \int_0^{\infty} \left\{ -\left(1 - e^{-x}\right)^n \right\} \frac{e^{-zx}}{x} dx$$

$$= \int_0^\infty \left\{ -(1-e^{-x})^{n+1} \right\} \frac{e^{-zx}}{(1-e^{-x})} \frac{dx}{x}$$

From here,

$$\sum_{n=1}^\infty \frac{1}{n+1} \sum_{k=0}^n \binom{n}{k} (-1)^k \log\left(1+\frac{k}{z}\right) = \int_0^\infty \left\{ -\sum_{n=1}^\infty \frac{(1-e^{-x})^{n+1}}{n+1} \right\} \frac{e^{-zx}}{(1-e^{-x})} \frac{dx}{x}$$

$$= \int_0^\infty \left\{ \log(1-(1-e^{-x})) + (1-e^{-x}) \right\} \frac{e^{-zx}}{(1-e^{-x})} \frac{dx}{x} = \int_0^\infty \left\{ -x + (1-e^{-x}) \right\} \frac{e^{-zx}}{(1-e^{-x})} \frac{dx}{x}$$

$$= \int_0^\infty \left\{ \frac{1}{x} - \frac{1}{1-e^{-x}} \right\} e^{-zx} dx = \psi(z) - \log z.$$

The summation in the very first sum can be started from $n=0$, as the first term there will be zero anyway. Thus

$$\psi(z) = \log z + \sum_{n=0}^\infty \frac{1}{n+1} \sum_{k=0}^n \binom{n}{k} (-1)^k \log\left(1+\frac{k}{z}\right). \tag{5.17}$$

Writing $\log\left(1+\frac{k}{z}\right) = \log(z+k) - \log z$ and noticing that for $n \geq 1$ we have

$$\sum_{k=0}^n \binom{n}{k} (-1)^k \log z = \log z \sum_{k=0}^n \binom{n}{k} (-1)^k = 0,$$

we conclude that

$$\sum_{n=0}^\infty \frac{1}{n+1} \sum_{k=0}^n \binom{n}{k} (-1)^k \log\left(1+\frac{k}{z}\right) = -\log z + \sum_{n=0}^\infty \frac{1}{n+1} \sum_{k=0}^n \binom{n}{k} (-1)^k \log(z+k)$$

and therefore,

$$\psi(z) = \sum_{n=0}^\infty \frac{1}{n+1} \sum_{k=0}^n \binom{n}{k} (-1)^k \log(z+k). \tag{5.18}$$

This representation can also be found in [9]. When $|t|<|z|$ the function $f(t)$ in (5.15) has the representation

$$f(t) = \log\left(1+\frac{t}{z}\right) = \sum_{m=1}^{\infty} \frac{(-1)^{m-1}}{m z^m} t^m,$$

and according to formula (2.24) we obtain from (5.17) the asymptotic expansion ($z \to \infty$)

$$\psi(z) = \log z + \sum_{m=1}^{\infty} \frac{(-1)^{m-1} B_m}{m} \frac{1}{z^m}.$$

For the same function $f(t)$ from (5.15) and with $y \geq 0$ we compute

$$\sum_{n=0}^{\infty} \frac{1}{n+1}\left\{\sum_{k=0}^{n}\binom{n}{k}(-1)^k f(y+k)\right\} = \sum_{n=0}^{\infty} \frac{1}{n+1}\left\{\sum_{k=0}^{n}\binom{n}{k}(-1)^k \log\left(1+\frac{y+k}{z}\right)\right\}$$

$$= \sum_{n=0}^{\infty} \frac{1}{n+1}\left\{\sum_{k=0}^{n}\binom{n}{k}(-1)^k \log(y+z+k)\right\} - \log z = \psi(y+z) - \log z,$$

where for the last equality we used (5.18). This leads, in view of (2.23), to the more general and interesting asymptotic representation ($z \to \infty$)

$$\psi(y+z) = \log z + \sum_{m=1}^{\infty} \frac{(-1)^{m-1} B_m(y)}{m} \frac{1}{z^m}. \tag{5.19}$$

As $B_1(y) = y - \frac{1}{2}$, we can write this in the form

$$\psi(y+z) = \log z + \frac{y}{z} - \frac{1}{2z} + \sum_{m=2}^{\infty} \frac{(-1)^{m-1} B_m(y)}{m} \frac{1}{z^m}. \tag{5.20}$$

Integrating for $z$ yields

$$\log \Gamma(y+z) = \left(z+y-\frac{1}{2}\right)\log z - z + \sum_{m=2}^{\infty} \frac{(-1)^{m-1} B_m(y)}{m(m-1)} \frac{1}{z^{m-1}} + C.$$

The constant of integration is $C = \log\sqrt{2\pi}$, as follows from Stirling's formula

$$\Gamma(z) \sim \sqrt{2\pi}\, z^{z-\frac{1}{2}} e^{-z} \quad (z \to \infty).$$

Finally,

$$\log \Gamma(y+z) = \left(z + y - \frac{1}{2}\right)\log z - z + \log\sqrt{2\pi} + \sum_{m=1}^{\infty} \frac{(-1)^m B_{m+1}(y)}{m(m+1)} \frac{1}{z^m}. \tag{5.21}$$

This is the classical asymptotic representation of the log-gamma function [14], [15], [18].

**Example 7** The Arakawa-Kaneko zeta function.

Recently in [1] and [8] the following function was introduced and studied

$$\zeta_r(s,a) = \frac{1}{\Gamma(s)} \int_0^\infty \mathrm{Li}_r(1-e^{-t}) \frac{t^{s-1} e^{-at}}{1-e^{-t}} dt, \tag{5.22}$$

where $a > 0$. In particular, it was proved in [1] and [8] that $\zeta_r(s,a)$ has the representation

$$\zeta_r(s,a) = \sum_{n=0}^{\infty} \frac{1}{(n+1)^r} \left\{ \sum_{k=0}^{n} \binom{n}{k} \frac{(-1)^k}{(k+a)^s} \right\}, \tag{5.23}$$

and when $r = 1$, it is related to the Hurwitz zeta function by the formula

$$\zeta_1(s,a) = s\zeta(s+1,a). \tag{5.24}$$

Applying equation (2.22) from Corollary 2 (with $z = 1$) to the function $f(t) = (t+a)^{-s}$ defined in (5.10) we find the representation

$$\zeta_r(s,a) = \sum_{m=0}^{\infty} \binom{-s}{m} \frac{(-1)^m B_m^{(r)}}{a^{m+s}} = \sum_{m=0}^{\infty} \binom{m+s-1}{m} \frac{B_m^{(r)}}{a^{m+s}}, \tag{5.25}$$

which is the asymptotic expansion of $\zeta_r(s,a)$ for $a$ large.

**References**

**[1]** **Bayad, Abdelmejid and Yoshinori Hamahata,** Polylogarithms and poly-Bernoulli polynomials, *Kyushu J. Math.*, 65 (2011) 15-24

**[2]** **Boyadzhiev, Khristo N.** Series Transformation Formulas of Euler Type, Hadamard Product of Functions, and Harmonic Number Identities**,** *Indian Journal of pure and Applied mathematics* 42 (2011) 371-387

**[3]** **Boyadzhiev, Khristo N.** A Series transformation formula and related polynomials, *Internat. J. Math. Math. Sci*.Vol. 2005 (2005), Issue 23, Pages 3849-3866

**[4]** **Boyadzhiev, Khristo N.** Evaluation of series with binomial sums, *Analysis Mathematica*, 40 (1), (2014), 13-23.

**[5]** **Boyadzhiev, Khristo N.** Exponential polynomials, Stirling numbers, and evaluation of some Gamma integrals, *Abstract and Applied Analysis*, Volume 2009, Article ID 168672 (electronic, open source).

**[6]** **Comtet**, **L.** Advanced Combinatorics (D. Reidel Publ. Co. Boston) (1974)

**[7]** **Connon, Donal F.** Some series and integrals involving the Riemann zeta function, binomial coefficients and the harmonic numbers. Volumes I and IV; http://arxiv.org/abs/0710.4022 and http://arxiv.org/abs/0710.4028

**[8]** **Coppo, M.-A. and B. Candelpergher**, The Arakawa-Kaneko zeta function, *Ramanujan J.*,22 (2010) 153-162

**[9]** **Guillera, Jesus, Jonathan Sondow**, Double integrals and infinite products for some classical constants via analytic continuations of Lerch's transcendent, *Ramanujan J*. 16 (2008) 247-270

**[10]** **Graham, Ronald L., Donald E. Knuth, and Oren Patashnik**, Concrete Mathematics, (Addison-Wesley Publ. Co., New York) (1994)

**[11]** **Hasse, H.**, Ein Summierungsverfahren fur die Riemannische $\zeta$ -Reihe, *Math. Z,* 32 (1930) 458-464

**[12]** **Jordan, Charles,** Calculus of finite differences (Chelsea, New York) (1950) (First edition: Budapest 1939).

**[13]** **Kaneko, M.** Poly-Bernoulli numbers. *J. Théorie de Nombres* , 9 (1997) 221–228

**[14]** **Magnus, W. Oberhettingeer, F. and R.P. Soni**, Formulas and Theorems for the Special Functions in Mathematical Physics (Springer) (1966)

**[15]** **Olver, Frank W.J. et al**. NIST Handbook of Mathematical Functions (NIST) (2010)

**[16]** **Todorov, Pavel G.** On the theory of the Bernoulli polynomials and numbers, *J. Math. Anal. Appl.* 104 (1984) 309-350

**[17]** **Todorov, Pavel G.** Taylor expansions of analytic functions related to $(1+z)^x - 1$, *J. Math.Anal. Appl.* 132 (1988) 264-280

**[18]** **Whittaker, E.T. and G.N. Watson**, A Course of Modern Analysis (Cambridge) (1992)